\documentclass[12pt, reqno]{amsart}
\usepackage{amsmath, amstext, amsbsy, amssymb, amscd}

\newtheorem{theorem}{Theorem}[section]
\newtheorem{lemma}[theorem]{Lemma}
\newtheorem{proposition}[theorem]{Proposition}
\newtheorem{corollary}[theorem]{Corollary}
\theoremstyle{definition}
\newtheorem{definition}[theorem]{Definition}

\theoremstyle{remark}
\newtheorem{remark}[theorem]{Remark}
\newtheorem{assumption}[theorem]{Assumption}
\newtheorem{conjecture}[theorem]{Conjecture}

\numberwithin{equation}{section}

\newcommand{\C}{\mathbb C}

\newcommand{\Hilb}{{\rm Hilb}}

\newcommand{\la}{\lambda}

\newcommand{\Quot}{{\rm Quot}}

\newcommand{\Supp}{{\rm Supp}}

\newcommand{\w}{\tilde}
\newcommand{\W}{\widetilde}

\newcommand{\Xn}{X^{[n]}}
\newcommand{\Z}{\mathbb Z}

{\vskip-\lastskip\medskip
  \noindent
  {\em #1.}\enspace
  }%
{\qed\par\medskip
  }

\begin{document}
\title[Euler numbers of moduli spaces] 
{On the Euler numbers of certain moduli spaces of curves 
and points}

\author[Wei-Ping Li]{Wei-Ping Li$^1$}
\address{Department of Mathematics, HKUST, Clear Water Bay, 
Kowloon, Hong Kong} \email{mawpli@ust.hk}
\thanks{${}^1$Partially supported by the grant HKUST6114/02P}

\author[Zhenbo Qin]{Zhenbo Qin$^2$}
\address{Department of Mathematics, University of Missouri, 
Columbia, MO 65211, USA} \email{zq@math.missouri.edu}
\thanks{${}^2$Partially supported by an NSF grant}


\begin{abstract}
We determine the topological Euler number of
certain moduli space of $1$-dimensional closed subschemes
in a smooth projective variety which admits a Zariski-locally
trivial fibration with $1$-dimensional fibers.
The main approach is to use virtual Hodge polynomials
and torus actions. The results might shed some light on 
the corresponding Donaldson-Thomas invariants.
\end{abstract}

\maketitle
\date{}

\section{Introduction}

Recently, there have been surging interests in studying
the moduli spaces of $1$-dimensional closed subschemes
in a smooth projective variety. The motivation 
comes from Donaldson-Thomas theory and its interplay
with Gromov-Witten theory and Gopakumar-Vafa invariants
\cite{MNOP1, MNOP2, Kat, OP} (see the references
there for other papers). Donaldson-Thomas theory was 
introduced in \cite{DT, Tho} via integrals over 
the moduli spaces of semistable sheaves and
via the theory of virtual fundamental cycles. 
It was further developed by 
Maulik-Pandharipande \cite{MP} and Jun Li. 
The moduli spaces of $1$-dimensional closed subschemes
in a smooth projective variety can be naturally regarded 
as the moduli spaces of rank-$1$ stable sheaves over 
the variety.

In \cite{MNOP1, MNOP2}, several interesting conjectures
regarding the interplay among Donaldson-Thomas theory,
Gromov-Witten theory and Gopakumar-Vafa invariants
have been proposed for $3$-folds. In particular, 
it was conjectured (see the Conjecture~2 in both papers) 
that the {\it reduced partition function} (with a formal
variable $q$) for the Donaldson-Thomas invariants 
is a rational function of $q$, and is invariant under
the transformation $q \to 1/q$ when the $3$-fold is 
Calabi-Yau. It is possible that, at least in the case of
Calabi-Yau $3$-folds, the Donaldson-Thomas invariants
are closely related to the topological Euler numbers
of the corresponding moduli spaces of $1$-dimensional 
closed subschemes. Hence we propose the following analogue 
to the Conjecture~2 in \cite{MNOP1, MNOP2}.

\begin{conjecture}  \label{conj_euler}
Let $X$ be a smooth projective complex variety.
Let $\mathfrak I_n(X, \beta)$ be the moduli space 
of $1$-dimensional closed subschemes $Z$ of $X$ 
satisfying (\ref{In}), and let $\Xn = \mathfrak I_n(X, 0)$ 
be the Hilbert scheme of length-$n$ $0$-dimensional 
closed subschemes of $X$. Then the {\it reduced} partition 
function for the Euler numbers
$${\sum_n \chi \big (\mathfrak I_n(X, \beta)\big ) \, q^n
\over \sum_n \chi(\Xn) \, q^n}$$ 
is a rational function of $q$, and is invariant under
$q \to 1/q$ when $K_X = 0$.
\end{conjecture}

While Donaldson-Thomas invariants
are difficult to calculate in general, there exist 
many effective methods to compute Euler numbers.
In this paper, we verify Conjecture~\ref{conj_euler}
under certain assumptions. Specifically,
we assume that $X$ admits a Zariski-locally trivial 
fibration $\mu: X \to S$ where $S$ is smooth, 
the fibers are smooth irreducible curves 
of genus-$g$, and $\beta \in H_2(X; \Z)$ is the class of 
a fiber. An element in the moduli space 
$\mathfrak I_{(1-g)+n}(X, \beta)$ consists of a fiber of 
$\mu$ together with some points (possibly embedded in 
the same fiber). For simplicity, let
\begin{eqnarray}   \label{def_Mn'}
\mathfrak M_n = \mathfrak I_{(1-g)+n}(X, \beta).
\end{eqnarray}

\begin{theorem}  \label{thm_euler'}
Let $r = \dim (X) \ge 2$, and $P_r(n)$ (respectively,
$\W P_r(n)$) be the number of $r$-dimensional (respectively,
punctual $r$-dimensional) partitions of $n$. Then, 
\begin{eqnarray*} 
\sum_{n =0}^{+\infty} \chi(\mathfrak M_n) \, q^n
= \sum_{n =0}^{+\infty} \chi(\Xn) \, q^n 
  \cdot \chi(S) \cdot 
  \left ( {\sum_{n=0}^{+\infty} \W P_r(n) \, q^n 
  \over \sum_{n=0}^{+\infty} P_r(n) \, q^n}
  \right )^{2-2g}.
\end{eqnarray*}
\end{theorem}

We remark that the Euler number of the Hilbert scheme
$\Xn$ has been determined by G\" ottsche and Cheah 
\cite{Go1, Che}, and $\chi(\mathfrak M_1)$ can be 
computed from the Lemma~1 in \cite{Kat} where 
the structure of $\mathfrak M_1$ has been described. 
Also, we refer to Definition~\ref{r-part} 
for the notion of $r$-dimensional (respectively, 
punctual $r$-dimensional) partitions of $n$. 
Furthermore, we conjecture (see Conjecture~\ref{P_WP}) that 
\begin{eqnarray}    \label{P_WP.1'}
{\sum_{n=0}^{+\infty} \W P_r(n) \, q^n 
  \over \sum_{n=0}^{+\infty} P_r(n) \, q^n}
= {1 \over (1-q)^{r-2}}.
\end{eqnarray}
This formula is trivially true when $r = 2$. 
The case $r = 3$ is proved in Lemma~\ref{lemma_r=3}. 
Note that $K_X = 0$ forces $g = 1$ under 
our assumption about $\mu$.

\begin{corollary}  \label{partial_conj}
Conjecture~\ref{conj_euler} is true when $2 \le r \le 3$
or $K_X = 0$.
\end{corollary}

Our main idea to prove Theorem~\ref{thm_euler'} is to 
use virtual Hodge polynomials, Cheah's combinatorial 
arguments \cite{Che}, and torus actions. More precisely,
we decompose the moduli space $\mathfrak M_n$ into
a disjoint union of locally closed subsets,
and prove that there exist bijective morphisms
between these locally closed subsets and 
certain spaces constructed from the local model 
$\C^{r-1} \times C$ of the fibration $\mu: X \to S$, 
where $C$ denotes a fixed fiber of $\mu$.
The spaces constructed from the local model 
$\C^{r-1} \times C$ consist of Zariski-locally trivial
fibrations involving simpler objects.
Our proofs here are parallel to those in \cite{Go1}.
Similar results hold for the Hilbert scheme $\Xn$.
Since virtual Hodge polynomials are preserved under 
decompositions and bijective morphisms,
these preparations enable us to reduce the computation 
of the virtual Hodge polynomial $e(\mathfrak M_n; s, t)$ 
to those of
\begin{eqnarray*}   
\mathfrak M_{n, C_0}^{\C^{r-1} \times C}, \qquad
(\C^{r-1} \times C)_{C_0}^{[n]}
\end{eqnarray*}
(see Definition~\ref{def1} and Remark~\ref{rmk0} for 
the notations), where $C_0 = \{O\} \times C$ and 
$O$ is the origin of $\C^{r-1}$. 
Applying Cheah's combinatorial arguments, 
we show that the virtual Hodge 
polynomials of $\mathfrak M_{n, C_0}^{\C^{r-1} \times C}$
and $(\C^{r-1} \times C)_{C_0}^{[n]}$ can be further reduced 
to those of the corresponding {\it punctual} spaces:
\begin{eqnarray*}   
\mathfrak M_{n, L, O}^{\C^r}, \qquad \Hilb^n(\C^r, O)
\end{eqnarray*}
where $L$ denotes a coordinate line in $\C^r$. 
Using the well-known fact that $e(\cdot \,; 1, 1) = \chi(\cdot)$,
we obtain a formula for $\chi(\mathfrak M_n)$ in terms
of the Euler numbers of $\mathfrak M_{n, L, O}^{\C^r}$
and $\Hilb^n(\C^r, O)$ (see Proposition~\ref{prop_euler}).
Finally, we take suitable torus actions
on $\mathfrak M_{n, L, O}^{\C^r}$ to determine its Euler 
number in terms of punctual $r$-dimensional partitions 
of $n$ (the torus actions on the punctual Hilbert scheme
$\Hilb^n(\C^r, O)$ were studied in \cite{Che}). 
This allows us to prove Theorem~\ref{thm_euler'}.

It would be interesting to see whether (\ref{P_WP.1'})
is true for $r > 3$. In addition, many results and 
techniques in this paper can be generalized to handle 
other situations, e.g., when the class $\beta$ is 
replaced by a higher multiple of the class of 
a fiber of $\mu$. We plan to discuss these 
generalizations in a forthcoming paper.

The paper is organized as follows.
In Sect.~\ref{sect_hodge}, the basic properties of 
virtual Hodge polynomials are reviewed.
In Sect.~\ref{sect_moduli}, we identify $\mathfrak M_n$
with certain Grothendieck Quot-scheme and introduce
a natural decomposition of $\mathfrak M_n$.
In Sect.~\ref{sect_bij}, various bijective morphisms
are constructed. In Sect.~\ref{sect_red}, we reduce 
the computation to the local model $\C^{r-1} \times C$.
In Sect.~\ref{sect_red_punc}, we further reduce 
the computation from the local model $\C^{r-1} \times C$ 
to the punctual spaces $\mathfrak M_{n, L, O}^{\C^r}$ 
and $\Hilb^n(\C^r, O)$. In Sect.~\ref{sect_act}, 
we use torus actions to verify Theorem~\ref{thm_euler'}
and (\ref{P_WP.1'}) when $r = 3$.

\bigskip\noindent
{\bf Acknowledgments.} The authors are grateful to 
Professors Robert Friedman, Sheldon Katz and Jun Li for 
valuable helps and stimulating discussions.
The first author thanks the Mathematics Department
at the University of Missouri for the hospitality and 
the support during his visit in January of 2005.
Thanks also go to the referee for suggesting improvements.
\section{\bf Virtual Hodge polynomials and Euler numbers} 
\label{sect_hodge}

Danilov and Khovanskii \cite{DK} introduced
virtual Hodge polynomials for reduced complex schemes.
These polynomials can be viewed as a convenient tool for 
computing the Hodge numbers of smooth projective varieties 
by reducing to computing those of simpler varieties.
They can also be used to compute Euler numbers.
In this section, we recall the basic properties of
virtual Hodge polynomials.

First of all, let $Y$ be a reduced complex scheme 
(not necessarily projective, irreducible or smooth). 
Mixed Hodge structures are defined on the cohomology 
$H_c^k(Y, \mathbb Q)$ with compact support 
(see \cite{Del, DK}). The mixed Hodge structures coincide 
with the classical one if $Y$ is projective and smooth. 
For each pair of integers $(m, n)$, 
define the virtual Hodge number
$$e^{m, n}(Y) = \sum_k (-1)^k h^{m, n}(H_c^k(Y, \mathbb Q)).$$ 
Then the {\it virtual Hodge polynomial} of $Y$ is 
defined to be
\begin{eqnarray}   \label{e_red}
e(Y; s, t) = \sum_{m, n} e^{m, n}(Y) s^m t^n.
\end{eqnarray}

Next, for an arbitrary complex scheme $Y$, we put 
\begin{eqnarray}   \label{e}
e(Y; s, t) = e(Y_{\text{red}}; s, t)
\end{eqnarray}
following \cite{Che}. By (\ref{e}) and the results 
in \cite{DK, Ful, Che} for reduced complex schemes, we see that 
virtual Hodge polynomials satisfy the following properties:

\begin{enumerate}
\item[(i)] When $Y$ is projective and smooth, $e(Y; s, t)$ is 
the usual Hodge polynomial of $Y$. 
For a general complex scheme $Y$, we have
\begin{eqnarray}   \label{hodge_i}
e(Y; 1, 1)= \chi(Y)
\end{eqnarray}
where $\chi(Y)$ denotes the topological Euler number of $Y$.

\item[(ii)] If $\displaystyle{Y = \coprod_{i=1}^n Y_i}$ is 
a finite disjoint union of locally closed subsets, then
\begin{eqnarray}   \label{hodge_ii}
e(Y; s, t) = \sum_{i=1}^n e(Y_i; s, t).
\end{eqnarray}

\item[(iii)] If $f: Y \to Y'$ is a Zariski-locally trivial bundle 
with fiber $F$, then
\begin{eqnarray}   \label{hodge_iii}
e(Y; s, t) = e(Y'; s, t) \cdot e(F; s, t).
\end{eqnarray}

\item[(iv)] If $f: Y \to Y'$ is a bijective morphism, then 
\begin{eqnarray}   \label{hodge_iv}
e(Y; s, t) = e(Y'; s, t).
\end{eqnarray}
\end{enumerate}

By the Lemma 5.6 in \cite{Che} and the Theorem 4.1
in \cite{LY}, if $Y$ is a reduced complex scheme with 
a $\C^*$-action and if $Y^{\C^*}$ denotes the set
of fixed points, then 
\begin{eqnarray}   \label{C*}
\chi(Y) = \chi \big (Y^{\C^*} \big ).
\end{eqnarray}

\section{\bf The moduli spaces $\mathfrak I_n(X, \beta)$ 
and $\mathfrak M_n$} 
\label{sect_moduli}

Let $X$ be a smooth projective complex variety of 
dimension $r$. For a fixed class $\beta \in H_2(X; \Z)$
and a fixed integer $n$, following the definitions and
notations in \cite{MNOP1, MNOP2},
we define $\mathfrak I_n(X, \beta)$ to be the moduli space 
of $1$-dimensional closed subschemes $Z$ of $X$ 
satisfying the two conditions:
\begin{eqnarray}   \label{In}
\chi(\mathcal O_Z)=n, \qquad [Z] = \beta
\end{eqnarray}
where $[Z]$ is the class associated to the dimension-$1$ 
component (weighted by their intrinsic multiplicities) of $Z$.
The degree-$0$ moduli space $\mathfrak I_n(X, 0)$ is isomorphic to 
the Hilbert scheme $\Xn$ parametrizing length-$n$ 
$0$-dimensional closed subschemes of $X$. In general, 
when $\beta \ne 0$, the space $\mathfrak I_n(X, \beta)$ is only
part of the Hilbert scheme defined in terms of certain
degree-$1$ Hilbert polynomial (see \cite{Gro}).
By the Lemma~1 in \cite{MNOP2}, when $\dim (X) = 3$, 
the virtual dimension of $\mathfrak I_n(X, \beta)$ is
\begin{eqnarray}   \label{dim_In}
-(\beta \cdot K_X).
\end{eqnarray}

In the rest of the paper, we adopt the following basic assumptions.

\begin{assumption} \label{assumption}
We assume that $X$ admits a Zariski-locally trivial fibration
\begin{eqnarray}   \label{fibration}
\mu: X \to S
\end{eqnarray}
where $S$ is smooth, the fibers are smooth irreducible curves 
of genus-$g$, and $\beta \in H_2(X; \Z)$ is the class of a fiber.
Then, $\mathfrak I_{(1-g)}(X, \beta) \cong S$. For $n \ge 0$, let
\begin{eqnarray}   \label{def_Mn}
\mathfrak M_n := \mathfrak I_{(1-g)+n}(X, \beta).
\end{eqnarray}
\end{assumption}
 
Our goal is to determine, under Assumption~\ref{assumption}, 
the partition function for the Euler numbers of the moduli spaces 
$\mathfrak M_n = \mathfrak I_{(1-g)+n}(X, \beta), n \ge 0$:
\begin{eqnarray}   \label{goal}
\sum_{n =0}^{+\infty} \chi \big ( \mathfrak M_n \big ) \, q^n.
\end{eqnarray}

\begin{lemma}  \label{Mn_Quot}
Let $\Delta$ be the diagonal of $S \times S$, and 
$\Quot^n_{(\text{\rm Id}_S \times \mu)^*I_\Delta/S \times X/S}$ 
be the Grothendieck Quot-scheme with the constant polynomial $n$.
Then, there exists an isomorphism between $\mathfrak M_n$ and 
$\Quot^n_{(\text{\rm Id}_S \times \mu)^*I_\Delta/S \times X/S}$.
\end{lemma}
\begin{proof}
Note that every element in $\mathfrak M_n = 
\mathfrak I_{(1-g)+n}(X, \beta)$ is of the form:
\begin{eqnarray}   \label{element}
Z = \Xi + \Theta
\end{eqnarray}
where $\Xi \in X^{[n-n_0]}$ for some $n_0$ satisfying 
$0 \le n_0 \le n$, $\Supp(\Xi) \cap \Supp(\Theta) = \emptyset$,
and the dimension-$1$ component $\Theta$ is equal to some curve
$C \in \mathfrak M_0$ together with embedded points of 
length-$n_0$ (i.e., $I_\Theta \subset I_C$ and the quotient
$I_C/I_\Theta$ is supported at finitely many points in $C$
with $h^0(X, I_C/I_\Theta) = n_0$). So we have a surjection
$$I_C \to I_C/I_Z \to 0$$
where the quotient $I_C/I_Z$ is supported at finitely many points,
and has length $n$. 

It follows that the universal quotient over 
$\Quot^n_{(\text{Id}_S \times \mu)^*I_\Delta/S \times X/S}
\times X$ induces a bijective morphism $\phi_1:
\Quot^n_{(\text{Id}_S \times \mu)^*I_\Delta/S \times X/S}
\to \mathfrak M_n$.

On the other hand, let $\mathcal I_n$ be the universal ideal
sheaf over $\mathfrak M_n \times X$. Let $\mathcal I_n'$
be the saturation of $\mathcal I_n \subset 
\mathcal O_{\mathfrak M_n \times X}$ 
(see Definition~1.1.5 in \cite{HL}).
Then, $\mathcal I_n'$ is a flat family of ideal sheaves in 
$\mathfrak M_0 \cong S$, and fits in an exact sequence
\begin{eqnarray*}   
0 \to \mathcal I_n \to \mathcal I_n' \to 
\mathcal Q \to 0
\end{eqnarray*}
over $\mathfrak M_n \times X$. Now the flat family 
$\mathcal I_n'$ and the quotient 
$\mathcal I_n' \to \mathcal Q \to 0$ induces a morphism 
$\phi_2: \mathfrak M_n \to 
\Quot^n_{(\text{Id}_S \times \mu)^*I_\Delta/S \times X/S}$
which is inverse to $\phi_1$.
\end{proof}

In view of Lemma~\ref{Mn_Quot}, we will make no difference
between $\mathfrak M_n$ and the Quot-scheme 
$\Quot^n_{(\text{Id}_S \times \mu)^*I_\Delta/S \times X/S}$.
In particular, we have a natural morphism:
\begin{eqnarray}   \label{Mn_S}
\mathfrak M_n \to S.
\end{eqnarray}
Moreover, over $\mathfrak M_n \times X$, there exists 
a universal quotient
\begin{eqnarray}   \label{Mn1}
p^*I_\Delta \to \mathcal Q \to 0
\end{eqnarray}
where $p$ is the composition of the morphism $\mathfrak M_n \times X
\to S \times X$ induced from (\ref{Mn_S}) and the morphism
$\text{Id}_S \times \mu: S \times X \to S \times S$.

\begin{definition}  \label{def1}
Let $\beta$ and $g$ be from Assumption~\ref{assumption}.
Let $0 \le n_0 \le n$.
\begin{enumerate}
\item[(i)] We define $\mathfrak M_{n, n_0}$ to be the locally 
closed subset of $\mathfrak M_n$ consisting of all the elements 
$Z = \Xi + \Theta$ from (\ref{element}) such that
$h^0(X, I_C/I_\Theta) = n_0$.

\item[(ii)] Fix a fiber $C$ of $\mu$ and a point $x \in C$.
We define $\mathfrak M_{n, C}$ 
(respectively, $\mathfrak M_{n, C, x}$) to be the closed 
subset of $\mathfrak M_{n, n}$ consisting of all the elements 
$Z = \Xi + \Theta$ from (\ref{element}) such that 
$\Xi = \emptyset$, $I_\Theta \subset I_C$, and
$I_C/I_\Theta$ is supported at finitely many points
(respectively, supported at the point $x$).

\item[(iii)] Fix a fiber $C$ of $\mu$, and a point $x \in X$.
Define $X^{[n], n_0}_C$ to be the locally closed subset of $X^{[n]}$ 
consisting of all the elements $\Xi = \Xi_1 + \Xi_2 \in X^{[n]}$
such that $\Supp(\Xi_1) \cap C = \emptyset$, 
$\Supp(\Xi_2) \subset C$, and $\ell(\Xi_2) = n_0$.
Define $X^{[n]}_C = X^{[n], n}_C$, and define $X^{[n]}_x$ to 
be the closed subset of $X^{[n]}$ consisting of 
all $\Xi \in X^{[n]}$ such that $\Supp(\Xi) = \{x\}$
(i.e., $X^{[n]}_x$ is the punctual Hilbert scheme at $x$).
\end{enumerate}
\end{definition}

\begin{remark}  \label{rmk0}
To emphasis the dependence on $X$, we will also denote 
the notations $\mathfrak M_n$, $\mathfrak M_{n, n_0}$, 
$\mathfrak M_{n, C}, \ldots$ by $\mathfrak M_n^X$, 
$\mathfrak M_{n, n_0}^X$, $\mathfrak M_{n, C}^X, \ldots$ 
respectively. 
\end{remark}
\section{\bf Various bijective morphisms} 
\label{sect_bij}

We begin with an outline of this section.
Consider the variety $X$ as in Assumption \ref{assumption}. 
By Definition~\ref{def1}~(i), $\mathfrak M_n$ has 
a decomposition of locally closed subsets:
\begin{eqnarray}   \label{decomp}
\mathfrak M_n = \coprod_{n_0 = 0}^n \mathfrak M_{n, n_0}.
\end{eqnarray}
In order to use the method of virtual Hodge polynomials, 
we need that there is a Zariski-locally trivial fibration space 
mapping to $\mathfrak M_{n, n_0}$  bijectively.  For instance, 
we see from (\ref{Mn_S}) that $\mathfrak M_{n, n}$ admits 
a morphism to $S$. The fiber of this map over $s \in S$ is 
isomorphic to $\mathfrak M_{n, C}$ where $C=\mu^{-1}(s)$ is 
a fiber of $\mu$: $X\to S$. If we consider another fiber
$C^{\prime}=\mu^{-1}(s^{\prime})$, then $\mathfrak M_{n, C}$ 
is isomorphic to $\mathfrak M_{n,C^{\prime}}$. 
In view of (\ref{hodge_iii}) and (\ref{hodge_iv}), 
it suffices to prove that $\mathfrak M_{n, C}\times S$ 
admits a bijective morphism to $\mathfrak M_{n,n}$
locally over $S$. To achieve this, we first consider 
the local model where we take $X=\mathbb C^{m}\times C$ 
in Lemma~\ref{Cr}. For the general case $\mu\colon X\to S$, 
we choose an open subset $U$ of $S$ such that there exists 
an \' etale morphism $U\to \mathbb C^m$. 
Then we prove in Lemma~\ref{u} the existence of 
a bijective morphism $\mathfrak M_{n, C}^{\mathbb C^m\times C}
\times_{\mathbb C^m} U\to \mathfrak M_{n, n}^{U\times C}$. 
For other strata $\mathfrak M_{n, n_0}$, Lemma~\ref{lemma_Mnn0}
provides a kind of locally trivial fibration description.  

Now we fix some notations. Let $C$ denote a fixed fiber of 
the fibration $\mu\colon X\to S$. Let $m = \dim(X) -1 =r-1$, 
$O$ be the origin of $\mathbb C^m$, and $C_0=\{O\}\times C$.

\begin{lemma}  \label{Cr}
There exists a bijective morphism over $\C^m$:
$$\Psi: \mathfrak M_{n, C_0}^{\C^m \times C} \times \C^m \to
\mathfrak M_{n, n}^{\C^m \times C}.$$
\end{lemma}
\begin{proof}
Let $X_0 = \C^m \times C$. Over $\mathfrak M_{n, C_0}^{X_0} \times X_0$, 
there exists a universal quotient
\begin{eqnarray}   \label{univ_quot}
p_2^*I_{C_0} \to \mathcal Q_0 \to 0
\end{eqnarray}
where $p_2: \mathfrak M_{n, C_0}^{X_0} \times X_0 \to X_0$ is 
the second projection. Let 
$$\sigma: \C^m \times \C^m \to \C^m$$
be the subtraction: $\sigma(u, v) = u- v$. Let $\Sigma =
\text{Id}_{\mathfrak M_{n, C_0}^{X_0}} \times \sigma 
\times \text{Id}_C$:
$$\mathfrak M_{n, C_0}^{X_0} \times \C^m \times X_0
= \mathfrak M_{n, C_0}^{X_0} \times \C^m \times \C^m \times C
\to \mathfrak M_{n, C_0}^{X_0} \times \C^m \times C
= \mathfrak M_{n, C_0}^{X_0} \times X_0.$$
Then we obtain a commutation diagram of morphisms: 
\begin{eqnarray*}
\begin{array}{ccc}
\mathfrak M_{n, C_0}^{X_0} \times \C^m \times X_0
   &{\overset \Sigma \to}& \mathfrak M_{n, C_0}^{X_0} \times X_0 \\
\downarrow \pi&&\downarrow\\
\C^m \times \C^m &{\overset \sigma \to}& \C^m
\end{array}
\end{eqnarray*}
where the two vertical morphisms are the natural projections. 
We have
\begin{eqnarray*}
\Sigma^*\big (p_2^*I_{C_0} \big ) \to \Sigma^*\mathcal Q_0 \to 0
\end{eqnarray*}
over $\mathfrak M_{n, C_0}^{X_0} \times \C^m \times X_0$.
Let $\Delta_0$ be the diagonal of $\C^m \times \C^m$. Then,
we see that
$\Sigma^*\big (p_2^*I_{C_0} \big ) = I_{\pi^{-1}\Delta_0}$.
Therefore, over $\mathfrak M_{n, C_0}^{X_0} \times \C^m \times X_0$, 
we have 
\begin{eqnarray}  \label{univ_quot2}
I_{\pi^{-1}\Delta_0} \to \Sigma^*\mathcal Q_0 \to 0.
\end{eqnarray}
Let $u \in \C^m$. The restriction of (\ref{univ_quot2})
to $\mathfrak M_{n, C_0}^{X_0} \times \{u\} \times X_0 \cong
\mathfrak M_{n, C_0}^{X_0} \times X_0$ is
$$p_2^*I_{\{u \} \times C} = \Sigma_u^*(p_2^*I_{C_0}) \to 
  \Sigma_u^*\mathcal Q_0 \to 0$$
where $\Sigma_u$ is the automorphism of 
$\mathfrak M_{n, C_0}^{X_0} \times X_0 = \mathfrak M_{n, C_0}^{X_0} 
\times \C^m \times C$ induced by 
$$\sigma_u: \C^m \to \C^m$$
with $\sigma_u(v) = u - v$. By the universal property,
(\ref{univ_quot2}) induces a morphism:
$$\Psi: \mathfrak M_{n, C_0}^{X_0} \times \C^m \to 
\mathfrak M_{n, n}^{X_0}.$$
The morphism $\Psi$ is bijective since every $\Sigma_u$ is 
an automorphism.
\end{proof}

\begin{lemma}  \label{u}
Let $f: U \to \C^m$ be an \' etale morphism. 
Then there exists a bijective morphism
$\W \Psi_f: \mathfrak M_{n, n}^{\C^m \times C} \times_{\C^m} U
\to \mathfrak M_{n, n}^{U \times C}$ over $U$.
\end{lemma}
\begin{proof}
Let $X = U \times C$ and $X_0 = \C^m \times C$. Then there exists
a universal quotient
\begin{eqnarray}  \label{u.1}
(\pi_0)^*I_{\Delta_0} \to \mathcal Q \to 0
\end{eqnarray}
over $\mathfrak M_{n, n}^{X_0} \times X_0$, where $\Delta_0$ is 
the diagonal of $\C^m \times \C^m$ and 
$\pi_0$ is the composition:
$$\mathfrak M_{n, n}^{X_0} \times X_0 \to \C^m \times X_0 
= \C^m \times (\C^m \times C) \to \C^m \times \C^m.$$
The projection $\mathfrak M_{n, n}^{X_0} \times_{\C^m} U \to 
\mathfrak M_{n, n}^{X_0}$ and the morphism $f \times \text{Id}_C:
X \to X_0$ induces
\begin{eqnarray*}
F: \big ( \mathfrak M_{n, n}^{X_0} \times_{\C^m} U \big )
\times X \to \mathfrak M_{n, n}^{X_0} \times X_0
\end{eqnarray*}
which can also be regarded as the base change of 
$f\times f: U\times U\to \mathbb C^m\times \mathbb C^m$ by $\pi_0$:
\begin{eqnarray}    \label{li4}
\begin{array}{ccc}
\big ( \mathfrak M_{n, n}^{X_0} \times_{\C^m} U \big ) \times X
  &\longrightarrow&U \times U\\
\qquad \downarrow F&&\qquad \,\,\,\, \downarrow f \times f\\
\mathfrak M_{n, n}^{X_0} \times X_0&\overset{\pi_0}{\longrightarrow}
  &\,\, \C^m \times \C^m.\\
\end{array}
\end{eqnarray}
Pulling-back the surjection (\ref{u.1}) via $F$, 
we obtain the surjection 
\begin{eqnarray}  \label{u.2}
F^*(\pi_0)^*I_{\Delta_0} \to F^*\mathcal Q \to 0
\end{eqnarray}
over $\big ( \mathfrak M_{n, n}^{X_0} \times_{\C^m} U \big )
\times X$. Regard $U \times U$ as a scheme over $\C^m$ by using 
$$U \times U \,\, {\overset {p_1} \to} \,\, U \,\,
{\overset f \to} \,\, \C^m$$
where $p_1$ is the first projection of $U \times U$. Then,
\begin{eqnarray}  \label{u.3}
\big ( \mathfrak M_{n, n}^{X_0} \times_{\C^m} U \big )
\times X = \mathfrak M_{n, n}^{X_0} \times_{\C^m} (U \times U)
\times C.
\end{eqnarray}
Let $\Delta$ be the diagonal of $U \times U$, and $\pi$ be
the composition:
$$\big ( \mathfrak M_{n, n}^{X_0} \times_{\C^m} U \big )
\times X = \mathfrak M_{n, n}^{X_0} \times_{\C^m} (U \times U)
\times C \to U \times U.$$ 
Restricting (\ref{u.2}) to $\pi^{-1}(\Delta) \subset 
\big ( \mathfrak M_{n, n}^{X_0} \times_{\C^m} U \big )
\times X$, we obtain
\begin{eqnarray}  \label{u.4}
F^*(\pi_0)^*I_{\Delta_0}|_{\pi^{-1}(\Delta)} \to 
F^*\mathcal Q |_{\pi^{-1}(\Delta)} \to 0.
\end{eqnarray}
Note that $\pi^{-1}(\Delta) \subset (\pi_0 \circ F)^{-1}(\Delta_0)$.
In fact, since $f$ is \' etale, $(\pi_0 \circ F)^{-1}(\Delta_0)$ 
is the disjoint union of $\pi^{-1}(\Delta)$ and 
some other irreducible components. Hence
\begin{eqnarray*}
F^*(\pi_0)^*I_{\Delta_0}|_{\pi^{-1}(\Delta)}
= I_{(\pi_0 \circ F)^{-1}(\Delta_0)}|_{\pi^{-1}(\Delta)}
= I_{\pi^{-1}(\Delta)}|_{\pi^{-1}(\Delta)}
= \pi^*I_\Delta|_{\pi^{-1}(\Delta)}.
\end{eqnarray*}
Using (\ref{u.4}) and the surjection $\pi^*I_\Delta \to 
\pi^*I_\Delta|_{\pi^{-1}(\Delta)}$, we obtain
\begin{eqnarray}  \label{u.5}
\pi^*I_\Delta \to F^*\mathcal Q |_{\pi^{-1}(\Delta)} \to 0
\end{eqnarray}
over $\big ( \mathfrak M_{n, n}^{X_0} \times_{\C^m} U \big )
\times X$. One checks that $F^*\mathcal Q |_{\pi^{-1}(\Delta)}$
is flat over $\mathfrak M_{n, n}^{X_0} \times_{\C^m} U$ and
that the quotient (\ref{u.5}) induces a morphism:
$$\W \Psi_f: \mathfrak M_{n, n}^{X_0} \times_{\C^m} U \to
\mathfrak M_{n, n}^{X}$$
over $U$. Using the completions of the points in $U$,
we see that $\W \Psi_f$ is bijective.
\end{proof}

\begin{proposition}  \label{prop:u}
Let $O$ be the origin of $\C^m$ and $C_0 = \{O\} \times C$.
Let $f: U \to \C^m$ be an \' etale morphism. 
Then there exists a bijective morphism over $U$:
$$\Psi_f: \mathfrak M_{n, C_0}^{\C^m \times C} \times U
\to \mathfrak M_{n, n}^{U \times C}.$$
\end{proposition}
\begin{proof}
Follows from Lemmas~\ref{Cr} and \ref{u} by putting 
$\Psi_f = \W \Psi_f \circ (\Psi \times_{\C^m} \text{Id}_U)$.
\end{proof}

\begin{remark}  \label{rmk:u}
The bijective morphism $\Psi_f$ in Proposition~\ref{prop:u} is 
in fact an isomorphism. To see this, we use an analytic open 
covering of $U$ to show that $\Psi_f$ is an isomorphism locally
in analytic category. This coupled with the bijectivity implies 
that the morphism $\Psi_f$ is indeed an isomorphism.
\end{remark}

\begin{definition}  \label{def2}
For $0 \le n_0 \le n$, define $Z_{n, i}$ to be the locally closed 
subset of $X^{[n]} \times S$ consisting of all the pairs 
$(\Xi, s)$ such that $\Xi = \Xi_1 + \Xi_2$ with $\Supp(\Xi_1)
\cap \mu^{-1}(s) = \emptyset$, $\Supp(\Xi_2) \subset \mu^{-1}(s)$,
and $\ell(\Xi_1) = i$. Put
$$W_n = Z_{n, n}, \qquad T_n = Z_{n, 0}.$$
\end{definition}

Note that we have natural morphisms $Z_{n, i} \to X^{[n]}$
and $Z_{n, i} \to S$.

\begin{lemma}  \label{lemma_Mnn0}
Let $0 \le n_0 \le n$. Then there exists 
a bijective morphism:
$$W_{n-n_0} \times_S \mathfrak M_{n_0, n_0} \to 
\mathfrak M_{n, n_0}.$$
\end{lemma}
\begin{proof}
Let $\pi_1$ and $\pi_2$ be the two natural projections of
$W_{n-n_0} \times_S \mathfrak M_{n_0, n_0}$. 
We see from (\ref{Mn1}) that over $\mathfrak M_{n_0, n_0} \times X$, 
there exists a universal quotient
\begin{eqnarray*}
p^*I_\Delta \to \mathcal Q_1 \to 0.
\end{eqnarray*}
So over $\big (W_{n-n_0} \times_S \mathfrak M_{n_0, n_0} \big ) 
\times X$, we have a surjection:
\begin{eqnarray}   \label{claim.1}
(\pi_2 \times \text{Id}_X)^*p^*I_\Delta \to 
(\pi_2 \times \text{Id}_X)^*\mathcal Q_1 \to 0.
\end{eqnarray}
In addition, over $X^{[n-n_0]} \times X$, 
we have a universal quotient
$$\mathcal O_{X^{[n-n_0]} \times X} \to \mathcal Q_2 \to 0.$$
Hence over $\big (W_{n-n_0} \times_S \mathfrak M_{n_0, n_0} \big ) 
\times X$, we have another surjection:
\begin{eqnarray}   \label{claim.2}
\mathcal O_{(W_{n-n_0} \times_S \mathfrak M_{n_0, n_0}) 
\times X} \to \pi^*\mathcal Q_2 \to 0
\end{eqnarray}
where $\pi$ is the composition of $\pi_1 \times \text{Id}_X:
\big (W_{n-n_0} \times_S \mathfrak M_{n_0, n_0} \big ) \times X
\to W_{n-n_0} \times X$ and the natural morphism
$W_{n-n_0} \times X \to X^{[n-n_0]} \times X$. Note that
\begin{eqnarray*}
&\Supp \big ((\pi_2 \times \text{Id}_X)^*\mathcal Q_1 \big ) 
  \subset (\pi_2 \times \text{Id}_X)^{-1}p^{-1}(\Delta),& \\
&(\pi_2 \times \text{Id}_X)^{-1}p^{-1}(\Delta) \cap 
  \Supp \big ( \pi^*\mathcal Q_2 \big ) = \emptyset.&
\end{eqnarray*}
Hence $(\pi_2 \times \text{Id}_X)^*\mathcal Q_1 \oplus
\pi^*\mathcal Q_2$ is flat over $W_{n-n_0} \times_S 
\mathfrak M_{n_0, n_0}$. Moreover, combining the two surjections
(\ref{claim.1}) and (\ref{claim.2}), we obtain a surjection:
\begin{eqnarray}   \label{claim.3}
(\pi_2 \times \text{Id}_X)^*p^*I_\Delta \to 
(\pi_2 \times \text{Id}_X)^*\mathcal Q_1 \oplus
\pi^*\mathcal Q_2 \to 0.
\end{eqnarray}
This surjection induces a morphism $\psi: W_{n-n_0} \times_S 
\mathfrak M_{n_0, n_0} \to \mathfrak M_n$. One checks that
$\text{im}(\psi) = \mathfrak M_{n, n_0}$ and that $\psi:
W_{n-n_0} \times_S \mathfrak M_{n_0, n_0} \to \mathfrak M_{n,n_0}$
is injective.
\end{proof}

\begin{remark}  \label{rmk1}
A similar argument shows that there exists a bijective morphism:
$$W_{n-n_0} \times_S T_{n_0} \to Z_{n, n_0}.$$
\end{remark}
\section{\bf Reduction to the local model $\C^{r-1} \times C$} 
\label{sect_red}

In this section, using the results proved in the previous section, 
we reduce the computation of the virtual Hodge 
polynomial of $\mathfrak M_n$ to those of $\Xn$, 
$\mathfrak M_{n, C_0}^{\mathbb C^{r-1}\times C}$,
and $(\mathbb C^{r-1}\times C)^{[n]}_{C_0}$
where $C_0=\{O\}\times C$ and $O$ is the origin of $\mathbb C^{r-1}$.

\begin{lemma}  \label{lemma_red}
Let $O$ be the origin of $\C^{r-1}$ and $C_0 = \{O\} \times C$. Then, 
\begin{eqnarray}   \label{lemma_red.1}
\sum_{n =0}^{+\infty} e(\mathfrak M_n; s, t) \, q^n
= \sum_{n =0}^{+\infty} e(W_n; s, t) \, q^n \cdot 
\sum_{n =0}^{+\infty} 
e(\mathfrak M_{n, C_0}^{\C^{r-1} \times C}; s, t) \, q^n.
\end{eqnarray}
\end{lemma}
\begin{proof}
By (\ref{decomp}), Lemma~\ref{lemma_Mnn0}, (\ref{hodge_ii}) 
and (\ref{hodge_iv}), we obtain:
\begin{eqnarray}    \label{lemma_red.2}
e(\mathfrak M_n; s, t) 
= \sum_{n_0 = 0}^n e(\mathfrak M_{n,n_0}; s, t)
= \sum_{n_0 = 0}^n e(W_{n-n_0} \times_S \mathfrak M_{n_0, n_0}; s, t).
\end{eqnarray} 
Consider the commutative diagram for the fiber product
$W_{n-n_0} \times_S \mathfrak M_{n_0, n_0}$:
\begin{eqnarray}    \label{lemma_red.3}
\begin{array}{ccc}
W_{n-n_0} \times_S \mathfrak M_{n_0, n_0} &\longrightarrow&
   \mathfrak M_{n_0, n_0} \\
\downarrow \phi_1&&\downarrow \phi_2\\
W_{n-n_0} &\longrightarrow&S.\\
\end{array}
\end{eqnarray}
By the Proposition~I.3.24 in \cite{Mil}, there exist an open
affine cover $\{ U_i \}_i$ of $S$ and \' etale morphisms
$f_i: U_i \to \C^{r-1}$. By Proposition~\ref{prop:u}, 
we see that for each $i$, there exists a bijective morphism 
over the open affine subset $U_i$:
$$\Psi_{f_i}: \mathfrak M_{n_0, C_0}^{\C^{r-1} \times C} \times U_i
\to (\phi_2)^{-1}(U_i).$$
So there exist a decomposition $S = \displaystyle{\coprod_i} S_i$
of locally closed subsets $S_i$ and bijective morphisms
$\Psi_{S_i}: \mathfrak M_{n_0, C_0}^{\C^{r-1} \times C} \times S_i
\to (\phi_2)^{-1}(S_i)$. By (\ref{lemma_red.3}), there exist 
a decomposition 
$$W_{n-n_0} = \displaystyle{\coprod_i} W_{n-n_0, i}$$
of locally closed subsets $W_{n-n_0, i}$ and bijective morphisms
$$\Psi_{W_{n-n_0, i}}: \mathfrak M_{n_0, C_0}^{\C^{r-1} \times C} 
\times W_{n-n_0, i} \to (\phi_1)^{-1}(W_{n-n_0, i}).$$ 
Combining this with (\ref{hodge_ii}) and (\ref{hodge_iv}), 
we conclude that
\begin{eqnarray}   \label{lemma_red.4}
   e(W_{n-n_0} \times_S \mathfrak M_{n_0, n_0}; s, t)
&=&\sum_{i} e \big ( (\phi_1)^{-1}(W_{n-n_0, i}); s, t \big )
   \nonumber   \\
&=&\sum_{i} e \big ( \mathfrak M_{n_0, C_0}^{\C^{r-1} \times C} 
   \times W_{n-n_0, i}; s, t \big ) \nonumber   \\
&=&\sum_{i} e( W_{n-n_0, i}; s, t) \cdot 
   e(\mathfrak M_{n_0, C_0}^{\C^{r-1} \times C}; s, t)
   \nonumber   \\
&=&e(W_{n-n_0}; s, t) \cdot 
   e(\mathfrak M_{n_0, C_0}^{\C^{r-1} \times C}; s, t).
\end{eqnarray} 
Now (\ref{lemma_red.1}) follows immediately from (\ref{lemma_red.2}) 
and (\ref{lemma_red.4}).
\end{proof}

\begin{lemma}  \label{XnS}
Let $O$ be the origin of $\C^{r-1}$ and $C_0 = \{O\} \times C$. Then, 
\begin{eqnarray}   \label{XnS.1}
\sum_{n =0}^{+\infty} e(X^{[n]} \times S; s, t) q^n 
= \sum_{n =0}^{+\infty} e(W_n; s, t) q^n \cdot
\sum_{n =0}^{+\infty} e \big (
(\C^{r-1} \times C)_{C_0}^{[n]}; s, t \big ) q^n.
\end{eqnarray} 
\end{lemma}
\begin{proof}
By Remark~\ref{rmk1}, we have an analogue of (\ref{lemma_red.2}):
\begin{eqnarray}   \label{XnS.2}
e(X^{[n]} \times S; s, t) = \sum_{n_0 = 0}^n 
e(W_{n-n_0} \times_S T_{n_0}; s, t).
\end{eqnarray} 
Let $\la$ be a partition of $n_0$, denoted by $\la \vdash n_0$.
Express $\la$ as $\la = (\la_1, \ldots, \la_\ell)$ where 
$\la_1 \ge \ldots \ge \la_\ell$ and 
$\la_1 + \ldots + \la_\ell = n_0$. 
We define $T_{\la}$ to be the locally closed subset of $T_{n_0}$ 
consisting of all the pairs $(\Xi, s)$ such that 
$\Xi = \Xi_1 + \ldots + \Xi_\ell$ where 
$$\Supp(\Xi_i) = \{x_i\} \subset \mu^{-1}(s),$$ 
$\ell(\Xi_i) = \la_i$, and 
the points $x_1, \ldots, x_\ell$ are distinct. Then,
\begin{eqnarray}   \label{XnS.3}
e(X^{[n]} \times S; s, t) = \sum_{n_0 = 0}^n 
\sum_{\la \vdash n_0} e(W_{n-n_0} \times_S T_\la; s, t).
\end{eqnarray} 
Using the Lemma~2.1.4 in \cite{Go2}, we can prove that 
the natural morphism $T_\la \to S$ is a Zariski-locally trivial 
fibration with fibers isomorphic to $(\C^{r-1} \times C)_{C_0}^\la$.
Here $(\C^{r-1} \times C)_{C_0}^\la$ denotes the locally closed
subset of $(\C^{r-1} \times C)_{C_0}^{[n_0]}$ consisting of 
$$\Xi' = \Xi_1' + \ldots + \Xi_\ell'$$
where $\Supp(\Xi_i') = \{x_i'\} \subset C_0$, 
$\ell(\Xi_i') = \la_i$, and $x_1', \ldots, x_\ell'$ are distinct.
Hence
\begin{eqnarray*}
   e(X^{[n]} \times S; s, t) 
&=&\sum_{n_0 = 0}^n \sum_{\la \vdash n_0} e(W_{n-n_0}; s, t)
   \cdot e \big ((\C^{r-1} \times C)_{C_0}^\la; s, t \big )    \\
&=&\sum_{n_0 = 0}^n e(W_{n-n_0}; s, t) \cdot 
   \sum_{\la \vdash n_0} e \big (
   (\C^{r-1} \times C)_{C_0}^\la; s, t \big )   \\
&=&\sum_{n_0 = 0}^n e(W_{n-n_0}; s, t) \cdot 
   e \big ((\C^{r-1} \times C)_{C_0}^{[n_0]}; s, t \big ),
\end{eqnarray*}
where we used the fact that $(\C^{r-1} \times C)_{C_0}^{[n_0]}$ is 
the disjoint union of the locally closed subsets 
$(\C^{r-1} \times C)_{C_0}^\la, \la \vdash n_0$.
Now (\ref{XnS.1}) follows immediately.
\end{proof}

\begin{proposition}  \label{prop_red}
Let $O$ be the origin of $\C^{r-1}$ and $C_0 = \{O\} \times C$. Then, 
\begin{eqnarray*}   
\sum_{n =0}^{+\infty} e(\mathfrak M_n; s, t) \, q^n
= \sum_{n =0}^{+\infty} e(\Xn; s, t) \, q^n \cdot e(S; s, t) 
  \cdot {\sum_{n =0}^{+\infty} 
  e(\mathfrak M_{n, C_0}^{\C^{r-1} \times C}; s, t) \, q^n 
  \over \sum_{n =0}^{+\infty} e \big (
  (\C^{r-1} \times C)_{C_0}^{[n]}; s, t \big )}.
\end{eqnarray*}
\end{proposition}
\begin{proof}
The formula follows from Lemma~\ref{lemma_red} and 
Lemma~\ref{XnS}.
\end{proof}
\section{\bf Reduction to the punctual cases} 
\label{sect_red_punc}

 From Proposition~\ref{prop_red}, we see that it suffices to 
compute the virtual Hodge polynomials of 
$\mathfrak M_{n, C_0}^{\mathbb C^{r-1}\times C}$ and 
$(\mathbb C^{r-1}\times C)^{[n]}_{C_0}$. 
These spaces are similar to the Hilbert scheme $\Xn$ in the 
sense that they are all built up from the punctual cases. 
Cheah developed a method of computing virtual Hodge 
polynomials to deal with this kind of situation.
In order to apply the method to 
$\mathfrak M_{n, C_0}^{\mathbb C^{r-1}\times C}$ and 
$(\mathbb C^{r-1}\times C)^{[n]}_{C_0}$, we first recall Cheah's 
original approach in \cite{Che} for the case of $\Xn$.

Let $\Hilb^n(\C^r, O)$ be the punctual Hilbert scheme of $\C^r$
at the origin. Then there exist unique rational numbers
$H_{\ell, m, n}$ such that
\begin{eqnarray}   \label{power_hil1}
\sum_{n=0}^{+\infty} e(\Hilb^n(\C^r, O); s, t) q^n 
= \prod_{\ell=1}^{+\infty} \prod_{m, n = 0}^{+\infty}
\left ( {1 \over {1 - q^\ell s^m t^n}} \right )^{H_{\ell, m, n}}
\end{eqnarray}
as elements in $\mathbb Q[s, t][[q]]$. Define $\mathfrak h_r(q,s,t)
\in \mathbb Q[s, t][[q]]$ to be the power series:
\begin{eqnarray}   \label{power_hil2}
\mathfrak h_r(q, s,t)
= \sum_{\ell=1}^{+\infty} \left ( \sum_{m, n = 0}^{+\infty}
H_{\ell, m, n} s^m t^n \right ) q^\ell.
\end{eqnarray}
Then the main result proved in \cite{Che} states that 
\begin{eqnarray}   \label{thm_che}
\sum_{n=0}^{+\infty} e(\Xn; s, t) q^n 
= \text{\rm exp} \left ( \sum_{n=1}^{+\infty} {1 \over n}
e(X; s^n, t^n) \mathfrak h_r(q^n, s^n, t^n) \right ).
\end{eqnarray}
The key ingredients in Cheah's proof of (\ref{thm_che}) can be 
summarized as follows:
\begin{enumerate}
\item[{\bf (A)}] Each element $\Xi \in X^{[n]}$ can be uniquely
decomposed into $\Xi^{(1)} + \ldots + \Xi^{(\ell)}$ where
every $\Xi^{(i)} \in X^{[n_i]}$ is supported at 
a single point in $X$, $n_1 + \ldots + n_\ell = n$, 
and the supports of $\Xi^{(1)}, \ldots, \Xi^{(\ell)}$ 
are mutually distinct.

\item[{\bf (B)}] Every $X^{[n]}_x$ is 
isomorphic to $\Hilb^n(\C^r, O)$. Let $X^{[n]}_{(n)}$ be 
the closed subscheme of $\Xn$ consisting of all $\Xi \in \Xn$
such that $\Supp(\Xi)$ is a single point of $X$.
Then the natural morphism $X^{[n]}_{(n)} \to X$
sending $\Xi \in X^{[n]}_{(n)}$ to $\Supp(\Xi) \in X$
is Zariski-locally trivial with fibers isomorphic to 
$\Hilb^n(\C^r, O)$.

\item[{\bf (C)}] Using certain combinatorial arguments 
independent of $X$, one reduces the computation to 
the virtual Hodge polynomials of $X$ 
and $\Hilb^n(\C^r, O)$ which contribute to the terms 
$e(X; s^n, t^n)$ and $\mathfrak h_r(q^n, s^n, t^n)$ in
(\ref{thm_che}) respectively. 
\end{enumerate}

It follows that we can apply Cheah's arguments to 
the computations of 
$$\sum_{n =0}^{+\infty} e \big (
(\C^{r-1} \times C)_{C_0}^{[n]}; s, t \big ) \, q^n, 
\qquad
\sum_{n =0}^{+\infty} 
e(\mathfrak M_{n, C_0}^{\C^{r-1} \times C}; s, t) \, q^n$$
in a straightforward fashion. 
For $\displaystyle{\sum_{n =0}^{+\infty} 
e \big ((\C^{r-1} \times C)_{C_0}^{[n]}; s, t \big ) \, q^n}$, 
we have
\begin{enumerate}
\item[{\bf (A1)}] Each element $\Xi \in (\C^{r-1} \times C)_{C_0}^{[n]}$ 
can be uniquely decomposed into 
$$\Xi^{(1)} + \ldots + \Xi^{(\ell)}$$ 
where each $\Xi^{(i)} \in (\C^{r-1} \times C)_{C_0}^{[n_i]}$ is 
supported at a single point in $C_0$, 
$$n_1 + \ldots + n_\ell = n,$$
and the supports of $\Xi^{(1)}, \ldots, \Xi^{(\ell)}$ 
are mutually distinct.

\item[{\bf (B1)}] Every 
$(\C^{r-1} \times C)_x^{[n]}, x \in C_0$ 
is isomorphic to $\Hilb^n(\C^r, O)$. The natural morphism 
$(\C^{r-1} \times C)_{C_0}^{[n]} \to C_0$ sending 
$\Xi \in (\C^{r-1} \times C)_{C_0}^{[n]}$ to $\Supp(\Xi) \in C_0$
is Zariski-locally trivial with fibers isomorphic to 
$\Hilb^n(\C^r, O)$.

\item[{\bf (C1)}] The same combinatorial arguments from {\bf (C)}
reduces the computation to the virtual Hodge polynomials of $C_0$ 
and $\Hilb^n(\C^r, O)$. 
\end{enumerate}
Therefore, we conclude as in (\ref{thm_che}) the following formula:
\begin{eqnarray}   \label{XnC}
\sum_{n =0}^{+\infty} 
  e \big ((\C^{r-1} \times C)_{C_0}^{[n]}; s, t \big ) \, q^n
= \text{\rm exp} \left ( \sum_{n=1}^{+\infty} {1 \over n}
  e(C_0; s^n, t^n) \mathfrak h_r(q^n, s^n, t^n) \right ).
\end{eqnarray}

Next, for the computation of $\displaystyle{\sum_{n =0}^{+\infty} 
e(\mathfrak M_{n, C_0}^{\C^{r-1} \times C}; s, t) \, q^n}$, we have
\begin{enumerate}
\item[{\bf (A2)}] Let $\Theta \in \mathfrak M_{n, C_0}^{\C^{r-1} 
\times C}$. By the definition of $\mathfrak M_{n, C_0}^{\C^{r-1} 
\times C}$, the quotient $I_{C_0}/I_\Theta$ is supported at 
finitely many points in $C_0$. Put 
$$I_{C_0}/I_\Theta = Q_1 \oplus \cdots \oplus Q_\ell$$
where each $Q_i$ is supported at a single point in $C_0$,
and the supports of $Q_1, \ldots, Q_\ell$ are mutually distinct.
Let $f: I_{C_0} \to I_{C_0}/I_\Theta$ be the quotient map.
For $1 \le i \le \ell$,
define the subscheme $\Theta^{(i)}$ by putting 
$$I_{\Theta^{(i)}} = f^{-1}(Q_i).$$
Then $\Theta \in \mathfrak M_{n, C_0}^{\C^{r-1} \times C}$
gives rise to $\Theta^{(1)}, \ldots, \Theta^{(\ell)}$. 
It is clear that the process can be reversed. Hence 
$\Theta \in \mathfrak M_{n, C_0}^{\C^{r-1} \times C}$ can be 
formally written as 
$$\Theta = \Theta^{(1)} + \ldots + \Theta^{(\ell)}$$
in a unique way,
where $\Theta^{(i)} \in \mathfrak M_{n_i, C_0}^{\C^{r-1} 
\times C}$ for $1 \le i \le \ell$, $n_1 + \ldots + n_\ell = n$,
each quotient $I_{C_0}/I_{\Theta_i}$ is supported at 
a single point in $C_0$, and the supports of the quotients
$I_{C_0}/I_{\Theta_1}, \ldots, I_{C_0}/I_{\Theta_\ell}$ 
are mutually distinct.

\item[{\bf (B2)}] Let $x \in C_0$. Since $C$ is a smooth curve 
in $X$, we have an isomorphism 
\begin{eqnarray}   \label{B2}
\mathfrak M_{n, C_0, x}^{\C^{r-1} \times C} \cong 
\mathfrak M_{n, L, O}^{\C^r}
\end{eqnarray}
between the punctual moduli spaces, where $L$ is a coordinate 
line in $\C^r$, $O$ is the origin of $\C^r$,
and $\mathfrak M_{n, L, O}^{\C^r}$ parametrizes all the 
$1$-dimensional closed subschemes $\Theta$ of $\C^r$ 
such that $I_\Theta \subset I_L$, 
$\Supp\big ( I_L/I_\Theta\big ) = \{O\}$, and 
$$h^0(\C^r, I_L/I_\Theta) = n.$$
Let $\mathfrak M_{(n), C_0}^{\C^{r-1} \times C}$ be the subset
of $\mathfrak M_{n, C_0}^{\C^{r-1} \times C}$ consisting of
all $\Theta \in \mathfrak M_{n, C_0}^{\C^{r-1} \times C}$
such that $\Supp(I_{C_0}/I_\Theta)$ is a single point in $C_0$.
By the construction in \cite{Gro}, there is a natural morphism
from $\mathfrak M_{n}^{\C^{r-1} \times C}$ to the $n$-th
symmetric product $\text{Sym}^n(\C^{r-1} \times C)$.
Its restriction to $\mathfrak M_{(n), C_0}^{\C^{r-1} \times C}$
gives rise to a morphism 
$$\phi: \mathfrak M_{(n), C_0}^{\C^{r-1} \times C} \to C_0.$$
An argument similar to the proof of Proposition~\ref{prop:u}
shows that there exist a decomposition of locally closed
subsets
$$C_0 = \coprod_i C_{0, i},$$
and bijective morphisms over the locally closed
subsets $C_{0, i}$:
$$\Phi_i: \mathfrak M_{n, L, O}^{\C^r} \times C_{0, i} \to 
\phi^{-1}(C_{0, i}).$$ 

\item[{\bf (C2)}] The same combinatorial arguments from {\bf (C)}
reduces the computation to the virtual Hodge polynomials of $C_0$ 
and $\mathfrak M_{n, L, O}^{\C^r}$. 
\end{enumerate}
Hence once again, we conclude as in (\ref{thm_che}) the following:
\begin{eqnarray}   \label{LO}
\sum_{n=0}^{+\infty} 
e(\mathfrak M_{n, C_0, x}^{\C^{r-1} \times C}; s, t) q^n 
= \text{\rm exp} \left ( \sum_{n=1}^{+\infty} {1 \over n}
e(C_0; s^n, t^n) \mathfrak c_r(q^n, s^n, t^n) \right ),
\end{eqnarray}
where the power series $\mathfrak c_r(q, s,t)
\in \mathbb Q[s, t][[q]]$ is defined by
\begin{eqnarray}   \label{power_curve2}
\mathfrak c_r(q, s,t)
= \sum_{\ell=1}^{+\infty} \left ( \sum_{m, n = 0}^{+\infty}
C_{\ell, m, n} s^m t^n \right ) q^\ell,
\end{eqnarray}
and the rational numbers $C_{\ell, m, n}$ are the unique rational numbers
such that
\begin{eqnarray}   \label{power_curve1}
\sum_{n=0}^{+\infty} e \big (\mathfrak M_{n, L, O}^{\C^r}; s, t \big) q^n 
= \prod_{\ell=1}^{+\infty} \prod_{m, n = 0}^{+\infty}
\left ( {1 \over {1 - q^\ell s^m t^n}} \right )^{C_{\ell, m, n}}.
\end{eqnarray}

\begin{lemma}  \label{vir_hodge}
Let $O$ be the origin of $\C^{r-1}$ and $C_0 = \{O\} \times C$. 
Let $\mathfrak h_r(q, s, t)$
and $\mathfrak c_r(q, s, t)$ be from (\ref{power_hil2})
and (\ref{power_curve2}). Then, 
$\displaystyle{\sum_{n =0}^{+\infty} e(\mathfrak M_n; s, t) \, q^n}$ 
is equal to
\begin{eqnarray*}  
&&\qquad \qquad \sum_{n =0}^{+\infty} e(\Xn; s, t) \, q^n 
  \cdot e(S; s, t)    \\
&\cdot&\text{\rm exp} \left ( \sum_{n=1}^{+\infty} {1 \over n}
  e(C_0; s^n, t^n) \big [\mathfrak c_r(q^n, s^n, t^n) - 
  \mathfrak h_r(q^n, s^n, t^n) \big ] \right ).
\end{eqnarray*}
\end{lemma}
\begin{proof}
Follows immediately from Lemma~\ref{prop_red}, (\ref{XnC}) 
and (\ref{LO}).
\end{proof}

\begin{proposition}  \label{prop_euler}
Under Assumption~\ref{assumption}, we have
\begin{eqnarray*} 
\sum_{n =0}^{+\infty} \chi(\mathfrak M_n) \, q^n
= \sum_{n =0}^{+\infty} \chi(\Xn) \, q^n 
\cdot \chi(S) \cdot
\left ( {\sum_{n=0}^{+\infty} 
    \chi \big (\mathfrak M_{n, L, O}^{\C^r} \big ) \, q^n 
  \over 
  \sum_{n=0}^{+\infty} \chi \big (\Hilb^n(\C^r, O) \big ) \, q^n}
\right )^{2-2g}.
\end{eqnarray*}
\end{proposition}
\noindent
{\it Proof.}
Follows from Lemma~\ref{vir_hodge}, (\ref{hodge_i}) and 
the observation that 
\begin{equation}
\text{\rm exp} \left ( \sum_{n=1}^{+\infty} {1 \over n}
  \mathfrak c_r(q^n, 1, 1)\right )
= \sum_{n=0}^{+\infty} 
  \chi \big (\mathfrak M_{n, L, O}^{\C^r} \big ) \, q^n.
\tag*{$\qed$}
\end{equation}

\section{\bf Torus actions on $\Hilb^n(\C^r, O)$ and 
$\mathfrak M_{n, L, O}^{\C^r}$} 
\label{sect_act}

The last step in our computation is to determine 
the Euler number of $\mathfrak M_{n, L, O}^{\C^r}$ 
(the Euler number of $\Hilb^n(\C^r, O)$ has been 
calculated in \cite{Che}). According to (\ref{C*}), 
we can make use of a suitable $\mathbb C^*$-action on 
$\mathbb C^r$ and compute the Euler number of the fixed
locus of the induced $\mathbb C^*$-action on 
$\mathfrak M_{n, L, O}^{\C^r}$.

The fixed points of torus actions on the spaces 
$\Hilb^n(\C^r, O)$ and $\mathfrak M_{n, L, O}^{\C^r}$ 
are closely related to multi-dimensional partitions of $n$.

\begin{definition}  \label{r-part}
(i) Let $r \ge 2$ and $n \ge 0$. {\it An $r$-dimensional 
partition} (respectively, {\it a punctual $r$-dimensional 
partition}) of $n$ is an array 
\begin{eqnarray}   \label{r-part.1}
(n_{i_1, \ldots, i_{r-1}})_{i_1, \ldots, i_{r-1}}
\end{eqnarray}
of nonnegative integers $n_{i_1, \ldots, i_{r-1}}$ 
indexed by the tuples 
\begin{eqnarray}   \label{r-part.2}
(i_1, \ldots, i_{r-1}) \in (\mathbb Z_{\ge 0})^{r-1}
\end{eqnarray}
\big (respectively, by the tuples $(i_1, \ldots, i_{r-1}) 
\in (\mathbb Z_{\ge 0})^{r-1} - \{(0, \ldots, 0)\} \big )$
such that 
\begin{eqnarray}   \label{r-part.3}
\sum_{i_1, \ldots, i_{r-1}} n_{i_1, \ldots, i_{r-1}} = n,
\end{eqnarray}
and $n_{i_1, \ldots, i_{r-1}} \ge n_{j_1, \ldots, j_{r-1}}$
whenever $i_1 \le j_1, \ldots, i_{r-1} \le j_{r-1}$.

(ii) We define $P_r(n)$ (respectively, $\W P_r(n)$) to be 
the number of $r$-dimensional (respectively, 
punctual $r$-dimensional) partitions of $n$.
\end{definition}

We remark that Definition~\ref{r-part.1}~(i) is consistent 
with the one used in \cite{MNOP1}, while our $r$-dimensional 
partitions are $(r-1)$-dimensional partitions in \cite{Che}.

Torus actions on the punctual Hilbert scheme $\Hilb^n(\C^r, O)$ 
have been studied in \cite{Che}.
Let $z_1, \ldots, z_r$ be the coordinate functions of $\C^r$.
Then $\C^*$ acts on $\C^r$ by
\begin{eqnarray}   \label{action}
t(z_1, \ldots, z_r) = (t^{w_1}z_1, \ldots, t^{w_r}z_r),
\qquad t \in \C^*.
\end{eqnarray}
This $\C^*$-action on $\C^r$ induces a $\C^*$-action on 
$\Hilb^n(\C^r, O)$. Now choose $w_1, \ldots, w_r \in \mathbb Z$ 
properly. Then the $\C^*$-fixed points in $\Hilb^n(\C^r, O)$
are precisely those corresponding to the colength-$n$ ideals of
$\C[z_1, \ldots, z_r]$ generated by monomials. 
These ideals are in one-to-one correspondence with 
$r$-dimensional partitions of $n$. Indeed, 
given an $r$-dimensional partition 
$(n_{i_1, \ldots, i_{r-1}})_{i_1, \ldots, i_{r-1} \ge 0}$ of $n$,
the ideal of $\C[z_1, \ldots, z_r]$ generated by the monomials
$z_1^{i_1} \cdots z_{r-1}^{i_{r-1}}z_r^{n_{i_1, \ldots, i_{r-1}}}$
has colength-$n$. Conversely, given a colength-$n$ ideal $I$ of
$\C[z_1, \ldots, z_r]$ generated by monomials,
we obtain an $r$-dimensional partition 
$(n_{i_1, \ldots, i_{r-1}})_{i_1, \ldots, i_{r-1} \ge 0}$ of $n$
by putting
\begin{eqnarray}   \label{I_to_part}
n_{i_1, \ldots, i_{r-1}} = \text{min}\{i_r|\, 
z_1^{i_1} \cdots z_{r-1}^{i_{r-1}}z_r^{i_r} \in I\}.
\end{eqnarray}
Therefore, by (\ref{C*}), we have (see the Proposition~5.1 in 
\cite{Che}):
\begin{eqnarray}   \label{euler_punc_hil}
\chi \big (\Hilb^n(\C^r, O) \big ) = P_r(n).
\end{eqnarray}

When $r = \dim (X) = 3$, torus actions on the moduli space 
$\mathfrak M_n$ for a toric variety $X$ have been studied 
in \cite{MNOP1}. For $r \ge 3$ and for torus actions on 
our punctual moduli space $\mathfrak M_{n, L, O}^{\C^r}$,
we choose the line $L$ to be defined by the equations:
$$z_1 = \ldots = z_{r-1} = 0.$$
Then the $\C^*$-action (\ref{action}) on $\C^r$ induces 
a $\C^*$-action on $\mathfrak M_{n, L, O}^{\C^r}$.
Again, choose the weights $w_1, \ldots, w_r \in \mathbb Z$ 
in (\ref{action}) properly. Then 
the $\C^*$-fixed points in $\mathfrak M_{n, L, O}^{\C^r}$
are precisely those corresponding to the ideals $I$ of
$\C[z_1, \ldots, z_r]$ such that $I$ is generated by monomials,
$I \subset (z_1, \ldots, z_{r-1})$, and
$$\dim_\C {(z_1, \ldots, z_{r-1}) \over I} = n.$$
As in the previous paragraph, we see that these ideals are 
in one-to-one correspondence with punctual $r$-dimensional 
partitions of $n$ (note that a linear basis of the ideal 
$(z_1, \ldots, z_{r-1})$ consists of all the monomials 
$z_1^{i_1} \cdots z_{r-1}^{i_{r-1}}z_r^{i_r}$ with
$$(i_1, \ldots, i_{r-1}) \in (\mathbb Z_{\ge 0})^{r-1} 
- \{(0, \ldots, 0)\},$$
and $i_r \ge 0$). Therefore, we obtain from (\ref{C*}) that 
\begin{eqnarray}   \label{euler_punc_curve}
\chi \big (\mathfrak M_{n, L, O}^{\C^r} \big ) = \W P_r(n).
\end{eqnarray}

\begin{theorem}  \label{thm_euler}
Under Assumption~\ref{assumption}, let $r \ge 2$. Then, 
\begin{eqnarray*} 
\sum_{n =0}^{+\infty} \chi(\mathfrak M_n) \, q^n
= \sum_{n =0}^{+\infty} \chi(\Xn) \, q^n 
  \cdot \chi(S) \cdot 
  \left ( {\sum_{n=0}^{+\infty} \W P_r(n) \, q^n 
  \over \sum_{n=0}^{+\infty} P_r(n) \, q^n}
  \right )^{2-2g}.
\end{eqnarray*}
\end{theorem}
\begin{proof}
The formula follows from Proposition~\ref{prop_euler},
(\ref{euler_punc_hil}) and (\ref{euler_punc_curve}).
\end{proof}

\begin{corollary}  \label{ellitpic}
If $\mu: X \to S$ is an elliptic fibration, then
\begin{equation}
\sum_{n =0}^{+\infty} \chi(\mathfrak M_n) \, q^n = 
\sum_{n =0}^{+\infty} \chi(\Xn) \, q^n \cdot \chi(S).  \tag*{$\qed$}
\end{equation}
\end{corollary}

\begin{conjecture}  \label{P_WP}
Let $r \ge 2$, and let $P_r(n)$ and $\W P_r(n)$ be defined
above. Then,
\begin{eqnarray}     \label{P_WP.1}
{\sum_{n=0}^{+\infty} \W P_r(n) \, q^n 
  \over \sum_{n=0}^{+\infty} P_r(n) \, q^n}
= {1 \over (1-q)^{r-2}}.
\end{eqnarray}
\end{conjecture}


The conjecture is clearly true when $r = 2$. 
The next lemma handles $r=3$.

\begin{lemma}  \label{lemma_r=3}
Conjecture~\ref{P_WP} holds when $r = 3$.
\end{lemma}
\begin{proof}
We shall use notations and results from Sect.~11.2 of \cite{And}.
Identify our $3$-dimensional partitions with 
the {\it plane partitions} there,
i.e., our $3$-dimensional partition
$(n_{i_1, i_2})_{i_1, i_2 \ge 0}$ is identified with 
the plane partition whose entry at the lattice point
$(i_1, i_2), i_1, i_2 \ge 0$ in the plane is equal to 
$n_{i_1, i_2}$. Similarly, our punctual $3$-dimensional 
partitions will correspond to 
the {\it punctual plane partitions}.

Let $\pi_\ell(n_1, \ldots, n_k; q)$
be the generating function for plane partitions 
with at most $\ell$ columns, at most $k$ rows,
and with $n_i$ being the first entry in the $i$-th row.
Then, 
\begin{eqnarray}   \label{lemma_r=3.1}
  \pi_{\ell+1}(n_1, \ldots, n_k; q)
= q^{n_1 + \ldots + n_k} \sum_{m_k = 0}^{n_k} 
  \sum_{m_{k-1}=m_k}^{n_{k-1}} \cdots
  \sum_{m_1=m_2}^{n_1} \pi_\ell(m_1, \ldots, m_k; q).
\end{eqnarray}
by the formula (11.2.1) in \cite{And}. Let $S_{k, \ell}(m, n)$ 
(respectively, $\W S_{k, \ell}(m, n)$) denote the set of 
plane partitions (respectively, punctual plane partitions)
of $m$ with at most $\ell$ columns, at most $k$ rows,
and with each entry $\le n$. 
Let $p_{k, \ell}(m, n) = |S_{k, \ell}(m, n)|$ and
$\w p_{k, \ell}(m, n)= |\W S_{k, \ell}(m, n)|$. 
Define two generating functions:
\begin{eqnarray*} 
    \pi_{k, \ell}(n; q) 
&=&\sum_{m=0}^{+\infty} p_{k, \ell}(m, n) \, q^m \\
    \w \pi_{k, \ell}(n; q) 
&=&\sum_{m=0}^{+\infty} \w p_{k, \ell}(m, n) \, q^m. 
\end{eqnarray*}
So $\pi_{+\infty, +\infty}(+\infty; q) = 
\displaystyle{\sum_{m=0}^{+\infty} P_3(m) \, q^m}$ and
$\w \pi_{+\infty, +\infty}(+\infty; q) = 
\displaystyle{\sum_{m=0}^{+\infty} \W P_3(m) \, q^m}$.
Define 
\begin{eqnarray}   \label{lemma_r=3.1.1}
(q)_i = (1-q)(1-q^2) \cdots (1-q^i)
\end{eqnarray}
for positive integers $i$. By the Theorem~11.2 in \cite{And}, 
\begin{eqnarray}   \label{lemma_r=3.2}
\pi_{k, \ell}(n; q) 
= {(q)_1(q)_2 \cdots (q)_{k-1} 
        \over (q)_\ell(q)_{\ell+1}\cdots (q)_{\ell+k-1}} 
  \cdot 
  {(q)_{n+\ell}(q)_{n+\ell+1} \cdots (q)_{n+\ell+k-1} 
        \over (q)_n(q)_{n+1}\cdots (q)_{n+k-1}}.
\end{eqnarray}

Let $\w \la \in \W S_{k, \ell}(m, n)$. By placing $n$
at the origin of the plane, we obtain $\la \in 
S_{k, \ell}(m+n, n)$. Conversely, if $\la \in 
S_{k, \ell}(m+n, n)$ and if the part of $\la$ at the origin 
is $n$, then by deleting the part at the origin, 
we obtain $\w \la \in \W S_{k, \ell}(m, n)$. Hence
\begin{eqnarray}   \label{lemma_r=3.3}
\w \pi_{k, \ell}(n; q)
= q^{-n} \sum_{m_k \le \ldots \le m_2 \le n}
  \pi_\ell(n, m_2, \ldots, m_k; q).
\end{eqnarray}

Setting $n_1 = \ldots = n_k = n$ in (\ref{lemma_r=3.1}), 
we conclude that
\begin{eqnarray*} 
\pi_{\ell+1}(\underbrace{n, \ldots, n}_{k \text{ copies}}; q) 
&=&q^{kn} \sum_{m_k \le \ldots \le m_2 \le m_1 \le n} 
   \pi_\ell(m_1, \ldots, m_k; q)   \\
&=&q^{kn} \sum_{m_1=0}^n \sum_{m_k \le \ldots \le m_2 \le m_1} 
   \pi_\ell(m_1, m_2, \ldots, m_k; q).
\end{eqnarray*}
Combining this with formula (\ref{lemma_r=3.3}), we see that
\begin{eqnarray}   \label{lemma_r=3.4}
\pi_{\ell+1}(\underbrace{n, \ldots, n}_{k \text{ copies}}; q)
= q^{kn} \sum_{m_1=0}^n q^{m_1} \cdot \w \pi_{k, \ell}(m_1; q).
\end{eqnarray}
On the other hand, by the formula (11.2.8) in \cite{And},
\begin{eqnarray*}   
\pi_{\ell+1}(\underbrace{n, \ldots, n}_{k \text{ copies}}; q)
= q^{kn} \cdot \pi_{k, \ell}(n; q).
\end{eqnarray*}
Therefore, we see immediately from (\ref{lemma_r=3.4}) that
\begin{eqnarray*}   
\pi_{k, \ell}(n; q) 
= \sum_{m_1=0}^n q^{m_1} \cdot \w \pi_{k, \ell}(m_1; q).
\end{eqnarray*}
Thus, $\pi_{k, \ell}(n; q) - \pi_{k, \ell}(n-1; q)
= q^n \cdot \w \pi_{k, \ell}(n; q)$, i.e., 
\begin{eqnarray}   \label{lemma_r=3.5}
\w \pi_{k, \ell}(n; q) 
= q^{-n} \left [ \pi_{k, \ell}(n; q) 
    - \pi_{k, \ell}(n-1; q) \right ].
\end{eqnarray}

To take the limits $k, \ell, n \to +\infty$, 
we assume $|q| < 1$ in the rest of the proof.
By (\ref{lemma_r=3.5}), (\ref{lemma_r=3.2}) and 
the definition of $(q)_i$ from (\ref{lemma_r=3.1.1}),
we have
\begin{eqnarray*}   
& &\w \pi_{k, \ell}(n; q)           \\
&=&\pi_{k, \ell}(n-1; q) \cdot {1 \over q^{n}} 
   \left \{ 
     {(1-q^{n+\ell})(1-q^{n+\ell+1}) \cdots (1-q^{n+\ell+k-1}) 
     \over 
     (1-q^{n})(1-q^{n+1}) \cdots (1-q^{n+k-1})} - 1 
   \right \}     \\
&=&\pi_{k, \ell}(n-1; q) \cdot 
   {1 + q + \cdots + q^{k-1} + O(q^n) + O(q^\ell) \over 
   (1-q^{n})(1-q^{n+1}) \cdots (1-q^{n+k-1})}.
\end{eqnarray*}
Now taking the limit $\ell, n \to +\infty$, we conclude that
\begin{eqnarray*}   
\w \pi_{k, +\infty}(+\infty; q)           
= \pi_{k, +\infty}(+\infty; q) \cdot 
   (1 + q + \cdots + q^{k-1}).
\end{eqnarray*}
Finally, letting $k \to +\infty$, we immediately obtain
\begin{eqnarray*}   
\w \pi_{+\infty, +\infty}(+\infty; q)           
= \pi_{+\infty, +\infty}(+\infty; q) \cdot 
   {1 \over 1-q},
\end{eqnarray*}
i.e., $\displaystyle{\sum_{m=0}^{+\infty} \W P_3(m) \, q^m}
= \displaystyle{\sum_{m=0}^{+\infty} P_3(m) \, q^m 
\cdot {1 \over 1-q}}$. This proves (\ref{P_WP.1}) for $r=3$.
\end{proof}

\begin{corollary}  \label{r=3}
Under Assumption~\ref{assumption}, let $2 \le r \le 3$. Then, 
\begin{eqnarray*} 
\sum_{n =0}^{+\infty} \chi(\mathfrak M_n) \, q^n
= \sum_{n =0}^{+\infty} \chi(\Xn) \, q^n \cdot \chi(S) 
\cdot {1 \over (1-q)^{(r-2) \cdot (2-2g)}}.
\end{eqnarray*}
\end{corollary}
\begin{proof}
Follows immediately from Theorem~\ref{thm_euler} and 
Lemma~\ref{lemma_r=3}.
\end{proof}

\end{document}